\documentclass{amsart}

\usepackage{amssymb}

\newcommand{\cO}{{\frak O}}
\newcommand{\disc}{\mathrm{disc}} 
\newcommand{\F}{\mathbf{F}_q} 
\newcommand{\gen}{{\mathrm{Gen} }}
\newcommand{\GL}{{\mathbf{GL} }}
\newcommand{\Pic}{{\mathrm{Pic}\, }}
\newcommand{\pp}{\mathfrak{p}} 
\newcommand{\SL}{{\mathbf{SL} }}
\newcommand{\SO}{{\mathbf{SO} }}
\newcommand{\Z}{\mathbf Z}
\renewcommand{\b}[1]{{\mathbf{#1}}}
\renewcommand{\H}{\mathbf{F}_{q^2}}
\renewcommand{\O}{{\mathbf{O} }}

\newtheorem{thm}{Theorem}[section]
\newtheorem{cor}[thm]{Corollary}
\newtheorem{prop}[thm]{Corollary}
\newtheorem{lem}[thm]{Lemma}

\theoremstyle{rem}

\theoremstyle{definition}

\numberwithin{equation}{section}
\numberwithin{thm}{section}

\begin{document}

\subjclass[2000]{Primary: 11E25, 11E12. Secondary: 11E41, 11D09}

\title{Representations of Definite Binary Quadratic Forms over $\F[t]$} 

\author{Jean Bureau}
\address{Mathematics Department, Louisiana State University, Baton Rouge, LA 70803-4918, U.S.A.}
\email{jbureau@math.lsu.edu} 
\author{Jorge Morales}
\address{Mathematics Department, Louisiana State University, Baton Rouge, LA 70803-4918, U.S.A.}
\email{morales@math.lsu.edu} 

\begin{abstract} In this paper, we prove that a binary definite
  quadratic form over $\F [t]$, where $q$ is odd, is completely
  determined up to equivalence by the polynomials it represents up to
  degree $3m-2$, where $m$ is the degree of its discriminant. We also
  characterize, when $q>13$, all the definite binary forms over $\F [t]$
  that have class number one. 
  \end{abstract}

\maketitle

\section{Introduction} \label{intro} It is a natural question to ask
whether binary definite quadratic forms over the polynomial ring $\F
[t]$ are determined, up to equivalence, by the set of polynomials they
represent. Here $\F$  is the finite field of order $q$ and $q$ is odd.

The analogous question over $\Z$ has been answered affirmatively -- with
the notable exception of the forms $X^2+3 Y^2$ and $X^2+XY+Y^2$, which
have the same representation set but are not equivalent -- by Watson
\cite{Watson:1979gf}. Several related results appear in the literature
as far back as the mid-nineteenth century (see \cite{Watson:1980ys}).

We begin with the easier question whether the discriminant
of a binary definite quadratic form over $\F [t]$ is determined by its
representation set. In the classical case over $\Z$, Schering
\cite{Schering:1859fk} showed that this is the case up to
powers of 2. The same type of ideas are used here to
show in the polynomial context that if $Q$ and
$Q'$ represent the same polynomials up to degree $3m-2$, where
$m=\max\{\deg \disc (Q),\deg \disc (Q')\}$, then $\disc (Q)=\disc (Q')$
(Proposition \ref{binary1}).

The main result
of this paper is that if $Q$ and $Q'$ have the same discriminant and
represent the same polynomials up to degree equal to their second
successive minimum, then they are equivalent (Theorem \ref{binaires}).
We show that if such forms were not equivalent, then there would be an
elliptic curve over $\F$  that has more rational points than allowed by
Hasse's bound. If the condition on the discriminants is omitted,
then having the same representation set up to degree $3m-2$ is enough to
conclude equivalence (Theorem \ref{main}).

The same questions can be asked for ternary definite quadratic forms. We
show that in this case, the representation {\em sets} (as opposed to the
representation {\em numbers}), are not enough in general to determine
the equivalence class. We do so by constructing a family of
counterexamples (Corollary \ref{counterexample}). It turns out, however,
that the representation {\em numbers}, that is the number of times that
each polynomial is represented, are sufficient to determine the
equivalence class of a ternary form, as it will be showed in an upcoming
paper \cite{Bureau:2007fk}.

Finally, in Section 6, we show, assuming $q >13$, that if a definite
binary quadratic form $Q$ has class number one (i.e. its genus contains
only one equivalence class), then $\deg \disc(Q)\le 2$ (Theorem
\ref{classnumber}).

We are indebted to the referee for her/his useful remarks.

\section{Notation and terminology}

The following notation will be in force throughout the paper:

\begin{itemize}
\item[$\F$  :] The finite field of order $q$. We always assume $q$ odd.
\item[$A$ :] The polynomial ring $\F[t]$
\item[$K$ :] The field of rational functions $\F(t)$
\item[$\delta$ :] A fixed non-square of $\F^\times$.
\end{itemize}

A quadratic form $Q$ over $A$ is a homogeneous polynomial
\[
Q=\sum_{1\le i,j\le n} m_{ij} X_iX_j,
\]
where $M=(m_{ij})$ is an $n\times n$ symmetric matrix with coefficients
in $A$. The group $\GL_n(A)$ acts by linear change of variables on the set
of such forms. Two forms in the same $\GL_n(A)$-orbit are called {\em
equivalent}. Two forms in the same $\SL_n(A)$-orbit are called {\em
properly equivalent}.

The {\em discriminant} of $Q$ is defined by
\[
\disc(Q)=(-1)^{n(n-1)/2}\det(M)
\]
as an element of $A/{\F^\times}^2$. This is an invariant
of the equivalence class of $Q$.

The {\em representation set} of $Q$ is the set of polynomials
\[
V(Q)=\{Q(\b x): \b x\in A^n\},
\]
and the {\em degree $k$ representation set} is
  \[
V_k(Q)=\{Q(\b x): \b x\in A^n,\,\deg Q(\b x)\le k\}.
\]

The form $Q$ is {\em definite} if it is anisotropic over
the field $K_\infty=\F((1/t))$. This implies in particular
that $n\le 4$. 

A definite quadratic form $Q$ is {\em reduced} if $\deg m_{ii}\le \deg
m_{jj}$ for $i\le j$ and $\deg m_{ij}< \deg m_{ii}$ for $i<j$. Gerstein
\cite{Gerstein:2003mz} showed that every definite quadratic form is
equivalent to a reduced form and that two reduced forms in the same
equivalence class differ at most by a transformation in $\GL_n(\F)$. In
particular, the increasing sequence of degrees of the diagonal terms of
a reduced form
\[
(\deg m_{11}, \deg m_{22},\ldots,\deg m_{nn})
\]
is an invariant of its equivalence class. This sequence
is called the {\em successive minima} of $Q$ and will be denoted
by $(\mu_1(Q),\mu_2(Q),\ldots,\mu_n(Q))$.

In the case of binary forms, which are the main topic of this paper,
we will often write
\[
Q=(a,b,c)
\]
for the quadratic form
\[
Q=aX^2 + 2b XY +c Y^2.
\]
For binary forms, it is easy to see that being definite means simply that
$\disc(Q)=b^2-ac$ has either odd degree or has even degree and non-square
leading coefficient. Also, $Q$ reduced translates into the condition
\begin{equation}\label{E:reduced}
\deg b <\deg a\le \deg c.
\end{equation}

If $Q=(a,b,c)$ is definite and reduced, then
\begin{equation}\label{E:deg} \deg Q(x,y)=\max\{2\deg x+\mu_1, 2\deg
y+\mu_2\} \end{equation} 
for all $x,y\in A$, where $\mu_1$ and $\mu_2$
are the successive minima. When $\mu_1$ and $\mu_2$ have distinct
parity, the equality \eqref{E:deg} follows immediately from
\eqref{E:reduced}. When $\mu_1$ and $\mu_2$ have the same parity,
\eqref{E:deg} follows from \eqref{E:reduced} together with the fact that
the leading coefficient of $-ac$ is a non-square by definiteness .

\section{Successive minima and discriminant} 

\begin{lem}\label{intermediate} Let $Q=(a,b,c)$ be a definite reduced form
with successive minima $\mu_1<\mu_2$. If $f\in A$ is represented by $Q$ and
$\mu_1\le \deg f <\mu_2$, then $f=r^2 a$ for some $r\in A$.
\end{lem}

\begin{proof} Write $
f=a r^2+ 2b rs +c s^2$, with $r,s\in A$. If $\deg f<\mu_2$, then by \eqref{E:deg}
we must have $s=0$, that is $f=r^2 a$.
\end{proof}

\begin{lem} \label{minimabinaires} Let $Q$ and $Q'$ be definite binary
forms over $A$ with discriminants $d$ and $d'$ respectively. Let
$m=\max\{\deg d,\deg d'\}$. If $V_m(Q)=V_m(Q')$, then
$\mu_i(Q')=\mu_i(Q)$ ($i=1,2$) and $\deg d=\deg d'$. Moreover, there are
reduced bases in which the diagonal entries of the matrices of $Q$ and
$Q'$ have the same leading coefficients.

\end{lem}
\begin{proof} Let $Q=(a,b,c)$ and $Q'=(a',b',c')$ be in
reduced form. Let $\mu_i=\mu_i(Q)$ and $\mu'_i=\mu'_i(Q)$ ($i=1,2)$.
  Since $a$ is represented by $Q'$, we clearly have $\mu'_1\leq \mu_1$. 
If $\mu'_2> \mu_2$, then 
\[
\mu'_1\le \mu_1\le \mu_2<\mu'_2,
\]
and applying Lemma \ref{intermediate} to $Q'$, we get
$a=a' r^2$ and $c=a' s^2$ for some $s,r\in A$. In particular,
$\mu_1\equiv \mu_2 \pmod 2$. Let $k=(\mu_2-\mu_1)/2$ and consider
the expression
\[
Q(t^kx,y)=t^{2k} a x^2+ 2 t^k b xy +c y^2
\]
with $x,y\in \F$. Using the inequality \eqref{E:reduced}, 
we see that the coefficient of degree $\mu_2$  of $Q(t^kx,y)$
is 
\begin{equation}\label{leading}
a_{\mu_1} x^2+ c_{\mu_2} y^2,
\end{equation}
where $a_{\mu_1}$ and $c_{\mu_2}$ are the leading coefficients 
of $a$ and $c$ respectively. Since $a_{\mu_1}c_{\mu_2}\neq 0$, the quadratic form \eqref{leading}
is non-degenerate over $\F$   and
therefore represents all elements of $\F^\times$. If we choose in particular
$x,y$ so that \eqref{leading} is not in the square class of
$a'_{\mu_1'}$, then $Q(t^kx,y)$ cannot be represented by $Q'$, since
otherwise it would be of the form $r^2a'$ by Lemma \ref{intermediate}.
Hence $\mu'_2\le \mu_2$, and by symmetry $\mu_1=\mu_1'$ and
$\mu_2=\mu_2'$. The equality $\deg d=\deg d'$ follows immediately.

We can assume without loss of generality that $a=a'$. It remains to see that
the leading coefficients of $c$ and $c'$ are in the same square class. 
When $\mu_1\equiv\mu_2\pmod 2$, the leading
coefficients of $c$ and of $c'$ are both in the square class of
$-\delta a_{\mu_1}$, where $\delta\in \F$ is a non-square. When
$\mu_1\not\equiv\mu_2\pmod 2$, the leading coefficient of any element in
$V(Q')$ whose degree has the same parity as $\mu_2$ must be in the same
square class as the leading coefficient of $c'$. This applies in
particular to $c$.

\end{proof}

\begin{lem}\label{bounded} Let $Q$ be a primitive definite binary
quadratic form over $A$ with discriminant $d$ and let $p$ be an
irreducible factor of $d$. Then $Q$ represents a polynomial not
divisible by $p$ of degree $< \deg d$. \end{lem} \begin{proof} Write
$Q$ in reduced form $Q=(a,b,c)$. Clearly either $a$ or $c$ satisfies the
condition.
\end{proof}

\begin{lem}\label{modp} Let $Q$ be a primitive definite binary
quadratic form over $A$ with discriminant $d$. Let $p\in A$.
Then each element of $V(Q)$ is congruent modulo $p$ to an element
in $V_{2\deg p+ \deg d -2}(Q)$.

\end{lem} \begin{proof} Let $\{e_1,e_2\}$ be a reduced basis for $Q$.
Each element of $V(Q)$ is congruent modulo $p$
to an element of the form $Q(x_1 e_1+x_2 e_2)$ with $\deg x_i\le \deg p -1$.
Clearly $\deg Q(x_1 e_1+x_2 e_2)\le 2(\deg p-1)+\mu_2(Q)\le 2(\deg p-1)+\deg d$.
\end{proof}

\begin{prop} \label{binary1} Let $Q$ and $Q'$ be definite 
binary quadratic form over $A$ with discriminants $d,d'$ respectively. Let
$m=\max\{\deg d,\deg d'\}$.  If $V_{3m-2}(Q)= V_{3m-2}(Q')$ then
$d'\in d\ \F^{\times 2}$. 
\end{prop}

\begin{proof} The statement is trivial if $m=0$, so we shall assume
through the proof that $m\ge 1$.

Notice that the equality of representation sets is preserved by scaling;
hence $Q$ and $Q'$ may be assumed primitive.

We shall prove that for each irreducible polynomial $ p\in A$:
\[
V_{3m-2}(Q)\subset V_{3m-2}(Q')\quad \mathrm{implies}\quad v_{
p}(d')\leq v_{ p}(d),
\]
where $v_{ p}(\cdot)$ denotes the $ p$-adic valuation. This will show
that $d=u d'$, where $u\in \F^\times$, and Lemma \ref{minimabinaires}
shows that $u$ must be a square.

Let $n=v_{ p}(d)$ and $n'=v_{ p}(d')$. If $\deg (p)>m$, then trivially
$n=n'=0$, so we may assume $\deg p \le m$.

Let $L$ be the $A$-lattice on which $Q$ is defined and let
$M=(p^nL^\sharp)\cap L$, where $L^\sharp$ is the dual lattice with
respect to $Q$. Then it is easy to see that the form $Q_0=p^{-n} Q|_M$
is integral and primitive and has discriminant $d$. By Lemma
\ref{bounded}, $Q_0$ represents a polynomial $u$ relatively prime to $p$
with $\deg u\le m-1$. It follows that $p^nu$ is represented by $Q$ and
since $\deg p^nu\le 2m-1\le 3m-2$ it must also be represented by $Q'$.
In particular, $p^nu$ must be represented $p$-adically by $Q'$. Over
$A_p$, the form $Q'$ is equivalent to a diagonal form $(a,0,p^{n'}b)$
where $a,b$ are $p$-adic units. Then there exist $x,y\in A_p$ such that
\begin{equation}\label{padic}
p^n u= a x^2+p^{n'}b y^2.
\end{equation}
It follows from \eqref{padic} that if $n'>n$, then $n=v_p(a x^2)\equiv
0\pmod 2$. Consider now the lattice $N=(p^{n/2}
L^\sharp)\cap L$ and let $Q_1=p^{-n}Q|_N$. One sees immediately that
$Q_1$ is primitive, integral and $\disc(Q_1)=p^{-n}d$, so $Q_1$ is
$p$-unimodular and thus $V(Q_1)$ contains representatives of all classes
modulo $p$. In particular, $Q_1$ represents a polynomial $w$ that is
relatively prime to $p$ and is in a different square class modulo $p$ as $a$.
Furthermore, by Lemma \ref{modp}, $w$ can be chosen so that $\deg w\le 2\deg p+ \deg
(p^{-n}d)-2$. 

The polynomial $f=p^n w$ is obviously represented by $Q$ and has
degree $\le 2\deg p+ \deg d -2\le 3m-2$, so it is also represented by
$Q'$. Writing $f$ as in \eqref{padic} and dividing by $p^n$ we see that
$w$ is in the same square class as $a$, which is a contradiction. Hence
$n'\le n$.

\end{proof}

\section{Forms with the same representation sets in small degree } 

\begin{thm} \label{binaires} Assume $q>3$. 
Let $Q$ and $Q'$ be two binary definite positive binary quadratic forms
over $A$ with the same discriminant and the same successive minima
sequence $(\mu_1,\mu_2)$. Suppose that $V_{\mu_2}(Q)= V_{\mu_2}(Q')$.
Then $Q$ and $Q'$ are equivalent.
\end{thm}

\begin{proof} Let $Q=(a,b,c)$ and $Q'=(a',b',c')$ be reduced forms.
There is no loss of generality in making the following assumptions:
$a=a'$ is monic and $c,c'$ have same leading coefficients. When
$\mu_1\equiv\mu_2\pmod2$, the leading coefficients of $c$ and $c'$ can
be assumed to be equal to $-\delta$, for the fixed non-square
$\delta\in \F$.\\

\textbf{1}. Suppose that $\mu_1\not\equiv\mu_2 \pmod 2$. Since
$c$ is also represented by $Q$, it is represented by $Q'$; hence,
there are $f \in A$ and $\beta\ in \F$ such that
$c=af^2+2b'f\beta +c'\beta ^2$. The different parity of the successive
minima implies that $\beta =\pm1$. By changing $b'$ into $-b'$ if
necessary, we can assume that $\beta =1$. Let $\varphi=\left(
\begin{matrix} 1&f\\ 0&1\\ \end{matrix} \right)\in \GL_2(\F).$
Then $Q'':=Q\circ\varphi=(a,b'',c')$, for some $b''\in A$.
Since $\det(\varphi)=1$, it follows that
$\disc(Q'')=\disc(Q)=\disc(Q')$; hence, $ac'-b''^2=ac'-b'^2$.
This leads to $b''=\pm b'$. \\

\textbf{2}. Suppose that $\mu_1\equiv\mu_2 \pmod2$ and that
$\mu_1 <\mu_2$. It follows from the equality of the discriminants 
that $\deg(c'-c)<\mathrm{max}\{\deg b,\deg b'\}<\deg a$.

If $b=b'=0$, we conclude immediately that $c=c'$ by the equality
of the discriminants. So we may assume $b\ne 0$. 

Consider all the elements $au^2+2bu+c\in V(Q)$ with $u\in \F$. By
assumption, the equation
\begin{equation} \label{equality1}
  au^2+2bu+c=ax^2+2b'xy+c'y^2
\end{equation}
   is always solvable for some $x=x_k t^k+x_{k-1} t^{k-1}+\cdots +x_0\in
A$, where $k=(\mu_2-\mu_1)/2$, and $y\in \F$.

Notice that for degree reasons, the polynomials $a$, $b$ and $c$ are
linearly independent over $\F$  (recall that we are assuming $b\ne 0$),
hence the left hand side of \eqref{equality1} takes exactly $q$ values
as $u$ runs over $\F$ . The equality of the leading coefficients in
\eqref{equality1} gives
\begin{equation}
\label{equality2} -\delta=x_k^2-\delta y^2. 
\end{equation}

It is a standard fact that the number of pairs $(x_k,y)$ satisfying 
\eqref{equality2} is $q+1$ (see e.g. \cite[Theorem 2.59]{Gerstein2008}). 
Notice that if $(x_k,y)$ is a solution of
\eqref{equality2}, then so is $(-x_k,y)$, thus the number of possible
$y$'s appearing in a solution of \eqref {equality2} is $(q-1)/2
+2=(q+3)/2$.

Since $q>(q+3)/2$ by hypothesis, there must be two different values of
$u$ on the left-hand side of \eqref{equality1} with the same $y$ on the
right-hand side. In other words, there exist $u,v \in \F$, $u\ne v$,
such that the system
\begin{equation} \label{aa} 
\left\{ \begin{matrix}
  au^2+2bu+c=ax^2+2b'xy+c'y^2\\
  av^2+2bv+c=az^2+2b'zy+c'y^2 
  \end{matrix} \right.
\end{equation}
has a solution $(x,y,z)$, with $x,z\in A$ and $y\in \F$.
By subtracting the two lines of \eqref{aa}, we get

\[
a(u^2-v^2)+2b(u-v)= a(x^2-z^2)+2b'(x-z)y
\]

  By degree considerations $x^2-z^2=u^2-v^2$ and hence $x$ and $z$
  are constant. In particular $x_k=0$ (since $k=(\mu_2-\mu_1)/2
  >0$) and hence, by \eqref{equality2}, we have $y^2=1$.
  
  Going back to \eqref{equality1}, we get
  \[
  a(u^2-x^2)+2(bu-b'xy)=c'-c
  \]
 
 As observed earlier,
$\deg(c'-c)<\mathrm{max}\{\deg b,\deg b'\}<\deg a$.
Thus the above equality implies $u^2=x^2$. Thus $2(bu-b'xy)=2u(b\pm
 b')=c'-c$. Replacing $b'$ by $-b'$ if necessary, we can assume
 $2u(b+b')=c'-c$. Multiplying by $b-b'$ gives $2u
 a(c-c')=2u(b^2-b'^2)=(c'-c)(b-b')$ by the equality of the
 discriminants. Degree considerations again imply $c=c'$ and
 $b=\pm b'$. \\

\textbf{3}. Suppose that $\mu_1= \mu_2=n$. Write 
\begin{align*}
  &a=t^n+a_{n-1}t^{n-1}\cdots+a_0\\ &c=-\delta t^n
  +c_{n-1}t^{n-1}+\cdots+c_0\\ &c'=-\delta t^n
  +c_{n-1}'t^{n-1}+\cdots+c'_0\\ &b=b_kt^k+\cdots+b_0\\
  &b'=b'_kt^k+\cdots+b'_0,\end{align*}
  where $k=\max\{\deg{b},\deg{b'}\}$. If $b=b'=0$,
  we are done, so we may assume $k\ge 0$ and $b'_k\ne 0$.
  Note that
since $\disc(Q)=\disc(Q')$, we have $\deg(c-c')<k$ as in the
previous case.

  Since $V_n(Q)= V_n(Q')$, for any pair $(u ,v )\in\F^2$, there exists a
  pair $(x ,y  )\in \F^2$ such that
\begin{equation}
\label{reps} Q(u ,v )=Q'(x ,y  ).
\end{equation}

 Taking the coefficients of $t^n$ and $t^k$ in the above polynomials, we
 get the system of quadrics:
\begin{equation}
\label{quadrics} \left\{ 
\begin{array}{l}
u^2-\delta v^2=x^2-\delta y^2\\
a_ku^2+2b_ku v  +c_kv^2=a_k x^2+2b'_k x y+c_k y^2,
\end{array}
\right. 
\end{equation}
which defines an algebraic curve $E$ in $\mathbf{P}^3$. For every
$(u,v)\in \F^2\setminus\{0\}$, there is $(x,y)\in \F^2\setminus\{0\}$
satisfying \eqref{quadrics}. Notice also that if a quadruplet
$(u,v,x,y)$ satisfies \eqref{quadrics}, so does $(u,v,-x,-y)$ and that
the two sides of the first equation are forms anisotropic over $\F$ , so
$|E(\F)|\ge 2(q+1)$.

If the curve $E$ given by \eqref{quadrics} were smooth, then it would be
an elliptic curve and by the Hasse estimate \cite[Ch.
V]{Silverman:1986lr} we would have $|E(\F)|\le 2\sqrt{q} + q+1$, which
would contradict the above count. Thus $E$ cannot be a smooth curve.

It is also known that the intersection of two quadric hypersurfaces, say
$Q_1=0, Q_2=0$, in $\mathbf{P}^m$ is a smooth variety of codimension $2$
if and only if the binary form $\det(X Q_1+Y Q_2)$ of degree $m+1$ has
no multiple factor (see e.g. \cite[Remark
1.13.1]{Colliot-Thelene:1987lr} or \cite[Chap. XIII \S
11]{Hodge:1952lr}). In the case of our system \eqref{quadrics}, by
computing explicitly the discriminant of $\det(X Q_1+Y Q_2)$, where
$Q_1$, $Q_2$ are the two quaternary quadratic forms of \eqref{quadrics},
we get the condition
\begin{equation}\label{E:smooth}
\delta^4(b_k-b'_k)^4(b_k+b'_k)^4\left((a_k\delta+c_k)^2-4\delta
b_k'^2\right)\left((a_k\delta+c_k)^2-4\delta b_k^2\right)=0.
\end{equation}

Since $\delta$ is not a square in $\F$  and $b_k'\ne 0$ by assumption, we
must have either $b_k=\pm b'_k$ or $b_k=0$ and $a_k\delta+c_k=0$. We
shall rule out the second possibility.

Since $V_n(Q)=
V_n(Q')$, these sets span the same $\F$ -subspace of $A$; in particular 
$b'$ must be an $\F$ -linear combination of $a$, $b$ and $c$. Write
\[
b'=\alpha a+\beta b+\gamma c,
\]
with $\alpha, \beta, \gamma \in \F$. 
Taking terms of degree $n$ gives
\[
0=\alpha-\delta \gamma,
\]
which implies 
 \[ 
 b'=\gamma (\delta a +c)+\beta b. 
 \]
Taking now terms of degree $k$ we get
\[ 
 b_k'=\gamma (\delta a_k +c_k) +\beta b_k. 
\]
If $b_k=0$ and $a_k\delta+c_k=0$, then $b_k'=0$, which is a
contradiction with our assumption.

Thus $b_k=\pm b'_k$ is the only possibility. Replacing $b$ by $-b$ if
needed, we shall assume $b_k=b'_k$.

We shall now show that $b=b'$. Suppose by contradiction that $b\ne b'$
and let $m=\deg(b-b')<k$. Then, by the equality of the discriminants,
$\deg(b^2-b'^2)=m+k=n+\deg(c-c')$ , which implies $\deg(c-c')<m$ and in
particular $c_m=c'_m$.

Exactly the same argument that showed $b_k^2={b'}_k^2$ (just replace $k$
by $m$ in \eqref{quadrics}) shows that $b_m^2=b_m'^2$. Now consider the
system

\begin{equation}
\label{quadrics2} \left\{ 
\begin{array}{l}
a_mu^2+2b_m u v  +c_m v^2=a_m x^2+2b'_m x y+c_m y^2\\
a_ku^2+2b_ku v  +c_kv^2=a_k x^2+2b'_k x y+c_k y^2.
\end{array}
\right. 
\end{equation}
Adding the two equations and combining the result with the first equation in \eqref{quadrics}
we get the system

\begin{equation}
\label{quadrics3} \left\{ 
\begin{array}{l}
u^2-\delta v^2=x^2-\delta y^2\\
(a_k+a_m)u^2+2(b_k+b_m) u v  +(c_k+c_m) v^2=\\
\qquad \qquad(a_k+a_m) x^2+2(b'_k+b'_m) x y +(c_k+c_m) y^2,
\end{array}
\right. 
\end{equation}
Applying one more time the rational-point counting argument, this time
to the above system, we conclude that $(b_m-b_k)^2=(b'_m-b'_k)^2$, which
yields $b_m b_k=b'_m b'_k$. Since $b_k=b_k'\ne 0$ we conclude
$b_m=b_m'$, which contradicts the hypothesis that $m=\deg(b-b')$. Hence
$b=b'$ as claimed.
\end{proof}

Finally, putting together Proposition \ref{binary1}, Lemma \ref{minimabinaires} and Theorem 
\ref{binaires}, we get our main result:

\begin{thm}\label{main}  Assume $q>3$. 
Let $Q$ and $Q'$ be definite binary quadratic forms over $A$ with
discriminants $d$ and $d'$ respectively. Let $m=\break \max\{\deg d,\deg d'\}$.
If $V_{3m-2}(Q)=V_{3m-2}(Q')$, then $Q$ and $Q'$ are equivalent.
\end{thm}

\section{The Ternary Case} \label{example}

In this section we give an example showing that in the case of ternary
definite forms over $A$, the representation {\em sets} in general do not
determine the discriminant, much less the equivalence class of the form
\footnote{ However, the representation {\em numbers} do determine the
equivalence class of such forms as showed in \cite{Bureau:2006fk},
\cite{Bureau:2007fk}.}.

\begin{lem} \label{repr} Let $Q_a=X^2+t Y^2-\delta(t+a^2)Z^2$, where
$a\in \F^\times$. Then a polynomial $f\in A$ is represented by $Q_a$
over $A$ if and only if it is represented by $Q_a$ over
$A_{(t)}=\F[[t]]$.
\end{lem}
\begin{proof} By \cite[Theorem 3.5]{Chan:2005ly}, the form $Q_a$ has
class number one, so a polynomial $f\in A$ is represented by $Q_a$ over
$A$ if and only if it is represented locally everywhere. At primes $\pp
$ not dividing $\disc(Q_a)=\delta t(t+a^2)$, $Q_a$ is unimodular and
isotropic, hence represents everything. At $\pp =(t+a^2)$, since
$t\equiv -a^2 \pmod{\pp}$, $Q_a$ is equivalent to
$X^2-Y^2-\delta(t+a^2)Z^2$ which also represents everything since
$X^2-Y^2$ already does so. Thus the only condition is at the prime $\pp
=(t)$ (the condition at $\infty$ is automatic by reciprocity).
\end{proof}

\begin{cor}\label{indrepr} For each $a\in \F^\times$, let $Q_a$ be
as in Lemma \ref{repr}. The representation set $V(Q_a)$ does not
depend upon the choice of $a$.
\end{cor} 
\begin{proof} 
By virtue of Lemma \ref{repr}, it is enough to notice that $Q_a$ is
equivalent to $X^2+tY^2-\delta Z^2$ over $\F[[t]]$,
which is independent of $a$.
\end{proof}

\begin{cor}\label{counterexample} Assume $q\ge 5$ and choose
$a,b\in \F^\times$ such that $a^2\ne b^2$. Then $V(Q_a)=V(Q_b)$ but
$\disc(Q_a)\ne \disc(Q_b)$.
\end{cor}
\begin{proof} Clear by Corollary \ref{indrepr}.
\end{proof}

\section{Primitive binary forms of class number one}

In this section, we characterize primitive binary quadratic forms
 over $A=\F [t]$ of class number one. Although it should
 be possible, in principle, 
 to deduce the results below from general formulas such as the ones in
 \cite{Korte:1984qv}, we prefer to give here a direct argument.

We begin by a statement
on orders in quadratic extensions of $K=\F (t)$.

\begin{prop}\label{P:classgroup} Let $D=f^2 D_0\in A$, where
$D_0$ is a square-free polynomial of either odd degree or of even degree
and non-square leading coefficient, and $f\in A$ is a monic polynomial.
Let $B=A[\sqrt{D}]$. Assume that $\Pic(B)$ is an abelian $2$-group and
has at most one cyclic component of order $4$ and all other components
of order $2$. Then

\begin{enumerate}
\item If $\deg D_0>0$ and $q>13$, then $D$ is square-free (i.e. $f=1$)
and $\deg D\le 2$.

\item If $\deg D_0=0$ and $q>5$, then $\deg D\le 2.$
\end{enumerate}
\end{prop}

\begin{proof}
Let $\cO=A[\sqrt{D_0}]$. Notice that $\cO$ is the maximal $A$-order
in the field $E=K(\sqrt{D_0})$ and that $f$ is the conductor of $B$ in $\cO$.

There is an exact sequence

\begin{equation}\label{E:picgen}
1\longrightarrow \frac{\cO^\times}{B^\times} \longrightarrow
\frac{(\cO/f\cO)^\times}{(A/f A)^\times} \longrightarrow \Pic(B)
\longrightarrow \Pic(\cO) \longrightarrow 1.
\end{equation}
\bigskip

{\bf 1.\/} Assume $\deg D_0>0$. Then $\cO^\times=B^\times=\F^\times$ and
we get a shorter exact sequence

\begin{equation}\label{E:pic}
1 \longrightarrow \frac{(\cO/f\cO)^\times}{(A/f A)^\times}
\longrightarrow \Pic(B) \longrightarrow \Pic(\cO) \longrightarrow 1.
\end{equation}

Let $h$ be the radical of $f$ (i.e. the product of all irreducible monic
divisors of $f$). The subgroup $(1+h\cO/f\cO)/(1+h A/f A)$ of
$(\cO/f\cO)^\times/(A/f A)^\times$ has order $q^{\deg f-\deg h}$ and is
a $2$-group by the exact sequence \eqref{E:pic}, so we must have $f=h$,
i.e. $f$ is square-free.

Let $\pi$ be an irreducible factor of $f$ of degree $d$. Then
$(\cO/\pi\cO)^\times/(A/\pi A)^\times$ is a direct factor of
$(\cO/f\cO)^\times/(A/f A)^\times$ and is cyclic of order $q^d-1$ or
$q^d+1$ (according to whether $\pi$ is split or inert in $E$) or is
isomorphic to the additive group $\mathbf{F}_{q^d}$ when $\pi$ is
ramified. Clearly the latter case is impossible since $q$ is odd and in
the first two cases we must have $q^d\pm 1=2$ or $4$, which is also
impossible when $q>5$. Hence $f=1$, $D$ is square-free and $B=\cO$.

Let $r$ be the number of irreducible factors of $D$. It is well-known
that the $2$-rank of $\Pic(\cO)$ is $r-1$. Hence, under our present
hypotheses, $|\Pic(\cO)|\le 2^r$. The order of $\Pic(\cO)$ is
essentially the class number $h_E$ of $E$; more precisely
$|\Pic(\cO)|=h_E$ if $\deg D$ is odd and $|\Pic(\cO)|=2h_E$ if $\deg D$
is even \cite[Proposition 14.7]{Rosen:2002qf}.

Using the lower bound for $h_E$ given by the Riemann Hypothesis
\cite[Proposition 5.11]{Rosen:2002qf}, we get

\[
\begin{aligned}
(\sqrt q -1)^{\deg D-1} && \le && 2^r &&\mathrm{if}&& \deg D && \mathrm{is\ odd};\\
(\sqrt q -1)^{\deg D-2} && \le && 2^{r-1} && \mathrm{if}&& \deg D && \mathrm{is\ even.}
\end{aligned}
\]
When $\deg D\ge 3$, using the above inequalities and the obvious fact
that $r\le \deg D$, we get easily the inequality $\log_2 (\sqrt q-1)\le
3/2$, which is impossible if $q>13$.

\bigskip
{\bf 2.\/} Assume $\deg D_0=0$ and $\deg f>0$. Then $\cO=\H[t]$, so
$\Pic(\cO)=\{1\}$, $\cO^\times=\H^\times $ and $B^\times=\F ^\times$.
The exact sequence \eqref{E:picgen} becomes

\begin{equation}\label{E:pic2}
1 \longrightarrow \frac{\H^\times}{\F^\times}\longrightarrow
\frac{(\cO/f\cO)^\times}{(A/f A)^\times} \longrightarrow \Pic(B)
\longrightarrow 1
\end{equation}

Let $p$ be the characteristic of $\F $. Taking $p$-parts in the sequence
above (i.e. tensoring by $\Z_p$), we get $[(\cO/f\cO)^\times/((A/f
A)^\times)]_p=0$. Exactly the same argument as in Case 1 shows that $f$
must be square-free. Hence
\begin{equation}\label{E:decomp}
\frac{(\cO/f\cO)^\times}{(A/f A)^\times}=\prod_{\pi | f}
\frac{(\cO/\pi\cO)^\times}{(A/\pi A)^\times},
\end{equation}
where $\pi$ runs over all irreducible monic divisors of $f$. 

Notice that the factors on the right-hand side of \eqref{E:decomp} are
cyclic of order $q^{\deg \pi}+1$ if $\deg \pi$ is odd, and $q^{\deg
\pi}-1$ if $\deg \pi$ is even.

Let $\pi$ be an irreducible factor of $f$ of even degree, say $\deg
\pi=2m$, then, by the exact sequence \eqref{E:pic2}, $(q^{2m}-1)/(q+1)$
must be a $2$-power $\le 4$. This is possible only when $m=1$ and $q=3$
or $q=5$. Similarly, if $\deg \pi$ is odd, say $\deg \pi=2m+1$, then
$(q^{2m+1}+1)/(q+1)$ must be a $2$-power, but it is always an odd
number, so the only possibility is $m=0$, i.e. $\deg \pi=1$. Thus, when
$q>5$, $f$ is a product of linear factors.

If $q+1$ is divisible by an odd prime $\ell$, then, since $\Pic(B)$ is a
$2$-group, taking $\ell$-parts in \eqref{E:pic2} shows that there must
be only one factor in the decomposition \eqref{E:decomp}, i.e. $f$ is
irreducible (necessarily linear as shown above).

The only case left is when $q+1$ is a $2$-power. Notice that the factors
on the right-hand side of \eqref{E:decomp} are all cyclic of order
$q+1$, since all the $\pi$'s are linear. By the hypothesis on $\Pic(B)$,
if there is more than one factor in \eqref{E:decomp}, then $q+1$ is a
$2$-power $\le 4$. This is impossible if $q>3$. Thus, also in this case,
$f$ has only one irreducible, necessarily linear, factor.
\end{proof}

Let $(V,Q)$ be a quadratic space over the field $K=\F (t)$. Let
$L\subset V$ be an $A$-lattice and let $\gen(L)$ be the set of lattices
of $V$ in the genus of $L$. The orthogonal group $\O(V,Q)$ acts on
$\gen(L,Q)$ and the number of orbits (which is well-known to be finite)
is called the {\em class number} of $L$ and will be denoted by $h(L,Q)$,
or simply $h(Q)$ when the underlying lattice is obvious. The number of
orbits of the action of the subgroup $\SO(V,Q)$ on $\gen(L,Q)$ will be
denoted by $h^+(L,Q)$. Since $\SO(V,Q)$ has index 2 in $\O(V,Q)$, we
have $h^+(L,Q)\le 2h(L,Q)$.

If $(L,Q)$ is primitive of rank 2, then $h^+(L,Q)$ depends only on
$D=\disc(L,Q)$. Indeed, let $G_D$ be the set of classes of primitive
binary quadratic forms of discriminant $D$ up to orientation-preserving
(i.e. of determinant 1) linear transformation. This set is a group for
Gaussian composition \cite{Kneser:1982zr} and there is a natural exact
sequence relating $G_D$ and $\Pic(B)$, where $B=A[\sqrt{D}]$, (see
\cite[\S 6]{Kneser:1982zr}), which in our situation is
\begin{equation}\label{E:comp}
1\longrightarrow \F ^\times/{\F ^\times}^2 \longrightarrow G_D
\longrightarrow \Pic(B)\to 1.
\end{equation}
The {\em principal genus} consists of forms in the genus of the norm
form $X^2-D Y^2$ of $B$, and their classes in $G_D$ form a subgroup
$G_D^0$. The different genera are cosets for this subgroup and hence
they have all the same number of classes, i.e. $h^+(L,Q)=|G_D^0|$ for
all primitive quadratic lattices  $(L,Q)$ of discriminant $D$. It is
also well-known (and easy to see) that $G_D/G_D^0$ is $2$-elementary.

\begin{thm}\label{classnumber} Let $Q$ be a definite primitive binary
quadratic form over $A$, where $q>13$. If $h(Q)=1$ then $\deg
\disc(Q)\le 2$.
\end{thm}

\begin{proof}
If $h(Q)=1$, then $h^+(Q)=|G_D^0|\le 2$ and by the remarks above $G_D$
is an abelian $2$-group with at most one cyclic component of order $4$
and all others of order $2$. So is $\Pic(B)$ by the exact sequence
\eqref{E:comp} and we conclude by Proposition \ref{P:classgroup}.
\end{proof}

\noindent \textbf{Remark.} Theorem \ref{classnumber} is incorrect
without the assumption $q>13$. Here is a counterexample for $q=13$.

Let $D=t(t^2-1)$ and let $E$ be the elliptic curve over
$\mathbf{F}_{13}$ given by the equation $s^2=D$. Let $B=A[\sqrt{D}]$.
Then $\Pic(B)=E(\mathbf{F}_{13})\cong \Z/2\Z\oplus\Z/4\Z$. It is easy to
see that the exact sequence \eqref{E:comp} is split in this case, so
$G_D^0=2G_D=2E(\mathbf{F}_{13})\cong\Z/2\Z$. Let $Q_0$ be a form whose
class $[Q_0]$ generates $G_D^0$. Then the genus of any form $Q$ of
discriminant $D$ consists of the classes $[Q]$ and $[Q']=[Q]+[Q_0]$ in
$G_D$. If $[Q]$ has order 4, then $[Q']=-[Q]$ i.e. $Q$ and $Q'$ are
(improperly) equivalent and thus $h(Q)=1$. An explicit example is
$Q=(t-5,4,-(t^2+5t+11))$, which corresponds to the point $P=(5,4)$ of
order 4 in $E(\mathbf{F}_{13})$.
\\

Theorem \ref{classnumber} gives a converse of a result of Chan-Daniels
\cite{Chan:2005ly}. We summarize this in the following statement:

\begin{prop} Assume $q>13$. A binary definite quadratic form $Q$
over $A$ of discriminant $D$ has class number one if and only if it
satisfies one of the following conditions:

\begin{enumerate}
\item $\deg D\le 1$.
\item $\deg D=2$ and $\mu_1(Q)=1$.
\item $\deg D=2$, $\mu_1(Q)=0$ and $D$ is reducible.
\end{enumerate}

\end{prop}
\begin{proof} The ``if'' part follows from \cite[Lemma 3.7]{Chan:2005ly}
and the ensuing remark. The ``only if'' part is a consequence of Theorem
\ref{classnumber}.
\end{proof}

\bibliographystyle{amsplain}
\providecommand{\bysame}{\leavevmode\hbox to3em{\hrulefill}\thinspace}
\providecommand{\MR}{\relax\ifhmode\unskip\space\fi MR }
\providecommand{\MRhref}[2]{%
  \href{http://www.ams.org/mathscinet-getitem?mr=#1}{#2}
}
\providecommand{\href}[2]{#2}

\end{document}